\documentclass[11pt, a4paper]{amsart}
\usepackage{a4wide}
\usepackage{amssymb}
\usepackage{amsthm}
\usepackage{amsfonts}
\newtheorem{lemma}{\bf LEMMA}[section]
\newtheorem{theorem}[lemma]{\bf THEOREM}
\newtheorem{corollary}[lemma]{\bf COROLLARY}

\newtheorem{fact}[lemma]{\bf FACT}

\newcommand{\restrict}{\mbox{$\mid$}}

\title{Order convergence and compactness}
\author{Dominic van der Zypen}
\begin{document}
\date{}  \maketitle \begin{quote}{\footnotesize {\bf R\'esum\'e.} Soit $(P,\leq)$
 un ensemble partiellement ordonn\'e et soit $\tau$
une topologie compacte sur $P$ qui est plus fine que la
topologie d'intervalles.
Alors $\tau$ est contenu dans la topologie de convergence d'ordre.
}\end{quote}

\parindent = 0mm
\parskip = 2 mm

\section{Topologies on a given poset}
On any given partially ordered set $(P,\leq)$ there are topologies
arising from the given order in a natural way (see also \cite{Compendium}).
Perhaps the best
known such topology is the {\it interval topology}. Set
${\mathcal{S}}^- =\{P\setminus (x]: x \in P\}$, and
${\mathcal{S}}^+ =\{P\setminus [x): x \in P\}$ where $(x]=\{y
\in P : y \leq x\}$ and $[x)=\{y \in P : y \geq x\}$. Then
${\mathcal{S}} ={\mathcal{S}}^- \cup {\mathcal{S}}^+$ is a
subbase for the interval topology $\tau_i(P)$ on $P$.

There is another natural way to endow an arbitrary poset $(P,\leq)$ with a topology. We
want to describe this topology in the following.

\def\F{{\mathcal F}}
\def\G{{\mathcal G}}
\def\U{{\mathcal U}}

A {\it (set) filter} $\F$ on $(P,\leq)$ is a nonempty subset of the powerset
of $P$ such that

$-$  $\emptyset \notin \F$

$-$ $U, V \in \F$ implies $U\cap V \in \F$

$-$ $U \in \F$ and $V \supseteq U$ imply $V \in \F$.

(Note that the above concept can of course be defined for arbitrary sets.) For
any subset $A\subseteq P$ let {\it the set of lower bounds of
A} be denoted by $A^l=\{x \in P : x \leq a \textrm{ for all } a \in A\}$ and
the set of upper bounds by $A^u=\{x \in P : x
\geq a \textrm{ for all } a \in A\}.$ If ${\mathcal S}$ is a collection of subsets of $P$ then we
set ${\mathcal S}^l=\bigcup \{S^l : S \in {\mathcal S}\}$, similarly set
${\mathcal S}^u=\bigcup \{S^u : S \in {\mathcal S}\}$.

Let $A\subseteq P$ be a subset of a poset $P$ and $y \in P$. We
say that $y$ is the infimum of $A$ if $y$ is the
greatest element of $A^l$ and write $\bigwedge A = y$. Dually we
define the supremum of $A$, written $\bigvee A$. Note that in
general, suprema and infima need not exist.

\def\orcon{{\dot \to}}
Let $\F$ be a filter on a poset $P$ and let $x\in P$. We say that
$\F$ {\bf order-converges to} $x$, in symbols $\F \orcon x$, if
$\bigwedge \F^u = x = \bigvee \F^l$. Note that the principal
ultrafilter consisting of the subsets of $P$ that contain $x$
order-converges to $x$.

Now we are able to define the {\bf order convergence topology} $\tau_o(P)$
(called order topology in \cite{Erne})  on any given
poset $P$ by:
\begin{quote}$\tau_o(P)=\{U \subseteq P : $ for any $x \in U$ and
any filter $\F$ with $\F \dot \to x$ we have $U \in \F \}.$
\end{quote}

It is straightforward to verify that this is a topology. Indeed,
$\tau_o(P)$ is the finest topology on $P$ such that order
convergence implies topological convergence (which is not hard to
prove either). We will make constant use of the following facts:

\begin{fact} \label{basicfacts}
Let $P$ be a poset, let $\F$ be a filter on $P$. Then:
\begin{enumerate}
\item $x \in \F^u \Leftrightarrow (x] \in \F$ and $x\in \F^l
\Leftrightarrow [x) \in \F.$
\item If $\F \orcon x$ then $\F^u \neq \emptyset \neq \F^l$.
\item Suppose $\F \orcon x$. If $x \not \leq a$ then $P\setminus (a] \in \F$. Dually
if $x \not \geq b$ then $P\setminus [b) \in \F$.
\item If $\F \orcon x$ and $\G$ is a filter on $P$ with $\G \supseteq \F$ then $\G \orcon x$.
\end{enumerate}
\end{fact}
\begin{proof} The proofs of assertions 1 and 2 are straightforward, and assertion 3
follows directly from \cite{Erne}, p.~3.
We prove assertion 4. Since $\G^u \supseteq \F^u$ it suffices to show that $\G^u
\subseteq [x)$ in order to get $\bigwedge \G^u = x$. Assume that there is $y \in \G^u \setminus[x)$.
By assertion 1, $(y] \in  \G$.  Since we have $x \not \leq y$,  we get
$P \setminus (y] \in \F \subseteq \G$  (by assertion 3). So
$(y] \cap (P \setminus (y])=\emptyset \in \G$, which is a contradiction. The statement $\bigvee
\G^l=x$ is proved similarly.
\end{proof}
\section{The result}
Note that \ref{basicfacts}, assertion 3 implies that for any poset $P$, the interval topology $\tau_i(P)$
is contained in the order convergence topology $\tau_o(P)$. The following theorem
connects the concepts of interval topology, order convergence and compactness.
\begin{theorem}\label{datheorem}
Let $(P,\leq)$ be a poset. If $\tau$ is a compact topology on $P$ such that
$\tau_i(P) \subseteq \tau$, then $\tau \subseteq \tau_o(P)$.
\end{theorem}
\begin{proof}
Suppose that $W \in \tau\setminus \tau_o(P)$. Then there is $x \in W$ and a filter $\F$ on $P$
such that $\F \orcon x$ and $W \notin \F$.

The strategy now is to find an ultrafilter on the closed set $Q:=P\setminus W$ of $(P,\tau)$ that does not
converge to any point of $Q$ with respect to the subspace topology of $(P,\tau)$ on $Q$.
This will imply that $Q$ is a non-compact
closed subset of $(P,\tau)$, which in turn implies that $(P,\tau)$ is
not compact.

Note that every element of $\F$
intersects $Q$ (otherwise we would have $W \in \F$). So $\F \cup \{Q\}$ is a filter base
which is contained in some ultrafilter
$\U$.  Moreover, by \ref{basicfacts}, assertion 4, the ultrafilter
$\U$ order-converges to $x$.

\def\UPW{\U\restrict_{Q}}
It is easy to check that $$\UPW=\{U \cap Q : U \in \U\}$$ is an ultrafilter on $Q$ (this
uses of course the fact that $Q$ is a member of $\U$).

\textbf{Claim:}  $\UPW$
does not converge to any $y \in Q$ with respect to $\tau\restrict_{Q}$, the topology
on $Q$ induced by $\tau$.

{\it Proof of Claim}: Pick any $y \in Q$. First, we know that $x\in W$ and $y \in Q$,
whence $x \neq y$. Suppose that the following holds in $P$: \begin{quote}
(A) \hspace*{0.5cm} For all $z\in \U^u$
we have $y \leq z$ and for all $z'\in \U^l$ we have $y \geq z'$. \end{quote}
Then by definition
of order convergence this would imply $y \leq x$, since $x = \bigwedge \U^u$, and similarly
we would get $y \geq x$, a contradiction to $x\neq y$. So, (A) must be false, and without loss
of generality we may assume that there is a $z_0 \in \U^u$ with $y \not \leq z_0$. By
\ref{basicfacts}, assertion 1, we get $(z_0] \in \U$ which implies
$$B:=(z_0] \cap Q \in \UPW.$$ Since $y \not \leq z_0$
we also have
$$y \in P \setminus (z_0]. \quad \; (\star)$$
Because $\tau$ contains the interval topology $\tau_i(P)$, statement $(\star)$
above implies that the set $$V:=(P\setminus (z_0])\cap Q
= Q\setminus B$$ is an open neighborhood of $y$ in $(Q, \tau\restrict_{Q})$. But since $B\in \UPW$
and $V=Q\setminus B$, we have
$V\notin \UPW$, so $\UPW$ does not converge to $y$
with respect to $\tau\restrict_{Q}$. Since $y \in Q$ was arbitrary, the claim is proved.

The claim now shows that $Q=P\setminus W$ is a closed, non-compact subset of
$(P,\tau)$. So $(P,\tau)$
cannot be compact.
\end{proof}
This theorem has a direct consequence for Priestley spaces, i.e. compact
totally order-disconnected ordered spaces as introduced in (\cite{H0}, \cite{H1}).
\begin{corollary}
If $(P,\tau,\leq)$ is a Priestley space, then $\tau \subseteq \tau_o(P)$.
\end{corollary}

{\sc Acknowledgement:} The author wishes to thank Jana Fla\v skov\'a for pointing out
a mistake in the proof of fact \ref{basicfacts} and providing a more direct argument.
{\footnotesize

\parskip=0mm

\vspace*{1cm}
D.~van der Zypen\\
Allianz Suisse Insurance Company\\
Bleicherweg 19, CH-8022 Zurich, Switzerland\\
{\tt dominic.zypen@gmail.com}

\end{document}